\documentclass[12pt,reqno]{amsart}
\renewcommand{\O}{\mathcal{O}}

\usepackage{amsmath}
\usepackage{graphicx}
\usepackage{amsfonts}
\usepackage{amsthm}
\usepackage{pictex}
\usepackage{a4wide}

\def\Ai{\mathop{\mathrm{Ai}}\nolimits}
\def\qed{\hspace*{\fill}$\Box$}

\newtheorem{theorem}{Theorem}[section]
\newtheorem{proposition}[theorem]{Proposition}

\newtheorem{corollary}[theorem]{Corollary}
\newtheorem{Definition}[theorem]{Definition}
\newenvironment{definition}{\begin{Definition}\rm}{\end{Definition}}
\newtheorem{Remark}[theorem]{Remark}

\numberwithin{equation}{section}
\newcommand{\C}{\mathbb{C}}

\renewcommand{\O}{\mathcal{O}}
\renewcommand{\Re}{{\rm Re} \,}

\newcommand{\supp}{{\rm supp} \,}
\newcommand{\diag}{{\rm diag} \,}
\newcommand{\ds}{\displaystyle}

\title[Type II Hermite-Pad\'e approximatin]{Type II Hermite-Pad\'e approximation
\\ to the exponential function}

\author{A.B.J. Kuijlaars}
\address{Department of Mathematics,
Katholieke Universiteit Leuven,
Celestijnenlaan 200 B,
B-3001 Leuven, BELGIUM}
\email{arno@wis.kuleuven.be}
\author{H. Stahl}
\address{FB II Mathematik,
Technische Fachhochschule Berlin,
Luxemburgerstrasse 10,
D-13353 Berlin, GERMANY}
\email{stahl@tfh-berlin.de}
\author{W. Van Assche}
\address{Department of Mathematics,
Katholieke Universiteit Leuven,
Celestijnenlaan 200 B,
B-3001 Leuven, BELGIUM}
\email{walter@wis.kuleuven.be}
\author{F. Wielonsky}
\address{Laboratoire de Math\'ematiques P.~Painlev\'e,
UMR CNRS 8524 - Bat.~M2,
Universit\'e des Sciences et Technologies Lille,
F-59655 Villeneuve d'Ascq Cedex FRANCE}
\email{Franck.Wielonsky@math.univ-lille1.fr}

\dedicatory{Dedicated to Nico Temme on the occasion of his 65th birthday.}
\date{\today}
\thanks{The authors are supported in part by INTAS Research Network NeCCA 03-51-6637.
A.B.J.K. and W.V.A. are also supported by FWO Research project G.0455.04,
by K.U.Leuven research grant OT/04/24, and by NATO Collaborative Linkage Grant
PST.CLG.979738. A.B.J.K is also supported by grant BFM2001-3878-C02-02 of the
Ministry of Science and Technology of Spain and by the European Science Foundation Program
Methods of Integrable Systems, Geometry, Applied Mathematics
(MISGAM) and the European Network in Geometry, Mathematical
Physics and Applications (ENIGMA)}

\begin{document}
\begin{abstract}
We obtain strong and uniform asymptotics in every domain of the complex
plane for the scaled polynomials $a (3nz)$, $b (3nz)$, and $c (3nz)$
where $a$, $b$, and $c$ are the type II Hermite-Pad\'e approximants to
the exponential function of respective degrees $2n+2$, $2n$ and $2n$,
defined by $a (z)e^{-z}-b (z)=\O (z^{3n+2})$ and
$a (z)e^{z}-c (z)=\O (z^{3n+2})$ as $z\to 0$. Our analysis relies on a
characterization of these polynomials in terms of a $3\times 3$
matrix Riemann-Hilbert problem which, as a consequence of the famous
Mahler relations, corresponds by a simple transformation to a similar
Riemann-Hilbert problem for type I Hermite-Pad\'e approximants. Due to
this relation, the study that was performed in previous work, based on
the Deift-Zhou steepest descent method for Riemann-Hilbert problems,
can be reused to establish our present results.
\end{abstract}
\maketitle

\section{Hermite-Pad\'e approximation}
In this paper we consider quadratic Hermite-Pad\'e approximation to the
exponential function. Type I quadratic Hermite-Pad\'e approximation to
the exponential function near $0$ consists of finding polynomials
$p_{n_1,n_2,n_3}$, $q_{n_1,n_2,n_3}$ and $r_{n_1,n_2,n_3}$
of degrees $n_1,n_2$ and $n_3$ respectively, such that
\[  p_{n_1,n_2,n_3}(z)e^{-z} + q_{n_1,n_2,n_3}(z) +
r_{n_1,n_2,n_3}(z)e^{z} = {\O}(z^{n_1+n_2+n_3+2}),
  \qquad z \to 0. \]
If we set the right hand side equal to zero and solve for $e^z$, then
we obtain the
algebraic function
\[   \frac{-q_{n_1,n_2,n_3}(z) \pm
\sqrt{q_{n_1,n_2,n_3}^2-4p_{n_1,n_2,n_3}(z)r_{n_1,n_2,n_3}(z)}}
    {2r_{n_1,n_2,n_3}(z)} \]
as an approximation to $e^z$. Type II Hermite-Pad\'e approximation is
simultaneous
rational approximation to $e^{-z}$ and $e^{z}$ and consists of finding
polynomials
$a_{n_1,n_2,n_3}, b_{n_1,n_2,n_3}$ and $c_{n_1,n_2,n_3}$ of degrees at most
$n_2+n_3+2$, $n_1+n_3$ and $n_1+n_2$ respectively, such that
\begin{equation} \label{typeIIcond}
    \begin{aligned}
     a_{n_1,n_2,n_3}(z) e^{-z} - b_{n_1,n_2,n_3}(z) &=
    \O(z^{n_1+n_2+n_3+2}), \qquad z \to 0, \\
     a_{n_1,n_2,n_3}(z) e^{z} - c_{n_1,n_2,n_3}(z) &=
    \O(z^{n_1+n_2+n_3+2}), \qquad z \to 0.
    \end{aligned}
\end{equation}
This gives the rational approximants
$b_{n_1,n_2,n_3}(z)/a_{n_1,n_2,n_3}(z)$ to $e^{-z}$ and
$c_{n_1,n_2,n_3}(z)/a_{n_1,n_2,n_3}(z)$ to $e^{z}$, and both rational
approximants have the
same denominator. It is well-known that for the case of exponentials,
all indices $n_{1}$, $n_{2}$, $n_{3}$
are normal, i.e., the polynomials $a_{n_1,n_2,n_3}$, $b_{n_1,n_2,n_3}$
and $c_{n_1,n_2,n_3}$ exist and are unique up to a normalization
constant with degrees exactly $n_2+n_3+2$, $n_1+n_3$ and $n_1+n_2$, see
\cite[Theorem 2.1, p. 129]{NS}.

Hermite-Pad\'e approximation to the exponential function have been of
interest since Hermite and have recently been investigated by
Borwein \cite{BOR}, Chudnovsky \cite{CHU}, and Wielonsky \cite{wiel1,wiel2}.
The asymptotic distribution of the zeros for the scaled type I
Hermite-Pad\'e polynomials,
\begin{equation} \label{typeIpolys}
    P_n(z) = p_{n,n,n}(3nz), \qquad Q_n(z) = q_{n,n,n}(3nz),
    \qquad R_n(z) = r_{n,n,n}(3nz)
\end{equation}
and their asymptotic behavior as $n \to \infty$, have recently been
studied in detail by Stahl \cite{stahl1,stahl2,stahl3,stahl4} and by
Kuijlaars, Van Assche, and Wielonsky \cite{KVW}, see also \cite{KSVW}.
In \cite{KVW} a Riemann-Hilbert problem for type I
Hermite-Pad\'e approximation was formulated. The asymptotic analysis of this
Riemann-Hilbert problem with the Deift-Zhou \cite{DZ} steepest
descent method  for oscillatory Riemann-Hilbert problems and Stahl's
geometric description of the problem, allowed the authors
of \cite{KVW} to find strong asymptotic formulas for the polynomials
(\ref{typeIpolys}) as well as for the type I remainder term
\begin{equation} \label{errortypeI}
    E_n(z) = P_n(z) e^{-3nz} + Q_n(z) + R_n(z) e^{3nz}
\end{equation}
that hold uniformly in every region of the complex plane.
The paper \cite{KVW} contained the first instance of a steepest descent
analysis for a $3 \times 3$ matrix valued Riemann-Hilbert problem.
It was followed by the papers \cite{ABK,BK2} which dealt the asymptotic
analysis of $3 \times 3$ matrix valued Riemann-Hilbert problems
arising in random matrix theory. A Riemann-Hilbert analysis for
rational interpolants for the exponential function was carried out
in \cite{wiel3}.

It is the aim of this paper to show that the Riemann-Hilbert analysis
of \cite{KVW} also produces the corresponding asymptotic results
for the scaled type II Hermite-Pad\'e polynomials
\begin{equation} \label{typeIIpolys}
    A_n(z) = a_{n,n,n}(3nz), \qquad B_n(z) = b_{n,n,n}(3nz),
    \qquad  C_n(z) = c_{n,n,n}(3nz),
\end{equation}
and for the type II remainder terms
\begin{equation} \label{error12}
    E^{(1)}_n(z) = A_n(z)e^{-3nz} -B_n(z),\qquad
    E^{(2)}_n(z) = A_n(z)e^{3nz} - C_n(z).
\end{equation}
This is due to the fact that the type II Hermite-Pad\'e polynomials
are characterized by a Riemann-Hilbert problem, which is directly related
to the Riemann-Hilbert problem for type I Hermite-Pad\'e polynomials.
This relation was first observed in \cite{VAGeKu} and also used in
\cite{BK1,DK}. See also Section 2 below.
So we follow the asymptotic analysis of \cite{KVW},
and the reader is advised to consult a copy of that paper too
when reading the proofs in this paper.

\begin{figure}[tb]
\begin{center}
\includegraphics[height=8cm,width=10cm]{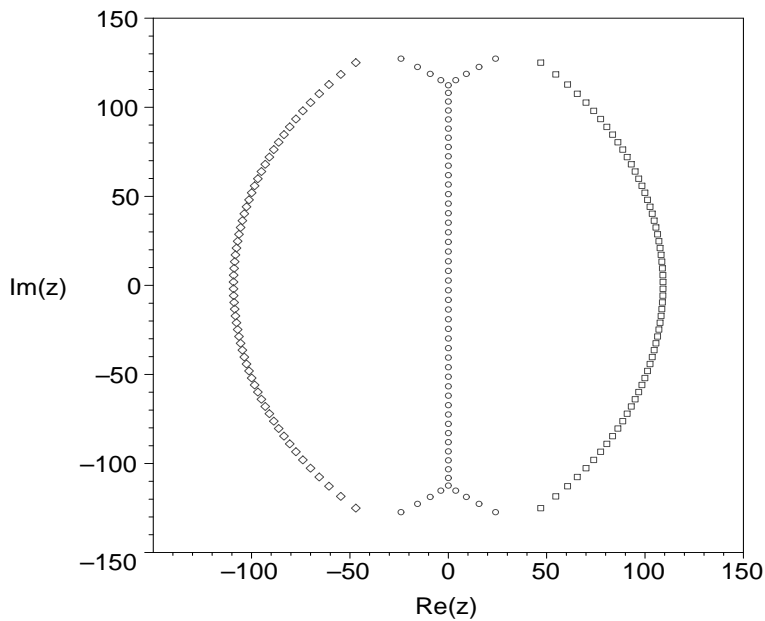}
\caption{Zeros of type I quadratic Hermite-Pad\'e polynomials
$p_{60,60,60}$ (the diamonds on the left) $q_{60,60,60}$ (the circles in
  the middle),
and $r_{60,60,60}$ (the boxes on the right)}
\label{zeroHPI}
\end{center}
\end{figure}
\begin{figure}[tb]
\begin{center}
\includegraphics[height=8cm,width=10cm]{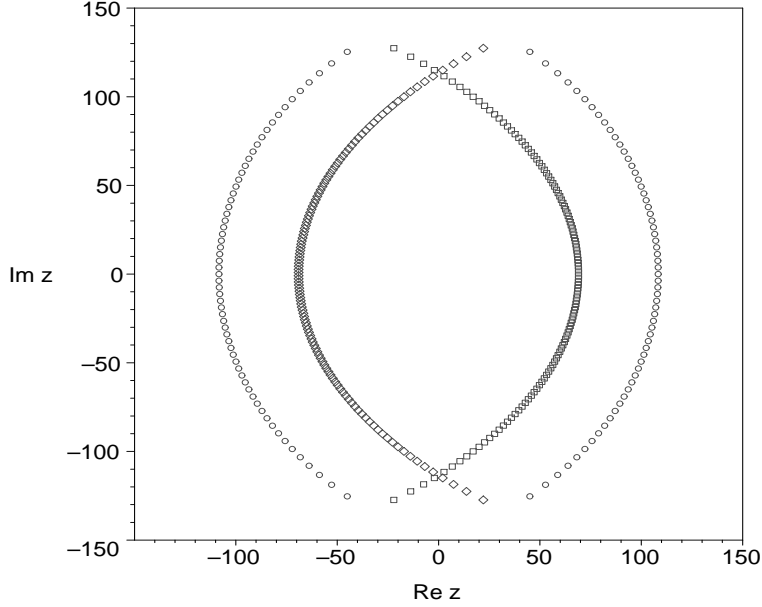}
\caption{Zeros of type II  quadratic Hermite--Pad\'e polynomials
$a_{60,60,60}$ (the circles on the left and right outer curves)
$b_{60,60,60}$ (boxes on the right inner curve),
and $c_{60,60,60}$ (diamonds on the left inner curve)}
\label{zeroHPII}
\end{center}
\end{figure}

To illustrate the connection between type I and type II we have
depicted in Figures \ref{zeroHPI} and \ref{zeroHPII}
the zeros of the type I and type II Hermite-Pad\'e approximants for
the case $n_1=n_2=n_3=60$. As it may be apparent from Figures 1
and 2, scaled zeros asymptotically accumulate on specific curves in the
complex plane, and the same system of curves is relevant
for both type I and type II approximants. More precisely,
the subsets of zeros on the left and on the right in Figures 1 and 2,
once scaled, tend to the same limit curves as the degree tends to
infinity. A similar assertion holds true for some of the zeros of
$q_{60,60,60}$ in the middle of Figure 1 (those outside of the imaginary
axis) and the corresponding subsets of zeros of $b_{60,60,60}$ and $c_{60,60,60}$
respectively lying in the left and right half-planes.

\section{Riemann-Hilbert problems for Hermite-Pad\'e approximation}
In \cite{KVW} we proved that the type I quadratic Hermite-Pad\'e
approximants for the exponential function may be characterized as the unique solution
of the following Riemann-Hilbert problem for a $3\times 3$ matrix valued function
$Y : \mathbb C \setminus \Gamma$, where $\Gamma$ is a closed contour
in the complex plane encircling the origin once in the positive direction,
cf. \cite[Theorem 5.1]{KVW},
\begin{enumerate}
\item[(1)] $Y$ is analytic in $\mathbb C \setminus \Gamma$
\item[(2)] $Y$ satisfies the jump condition
\begin{equation}  \label{eq:Yjump}
   Y_+(z) = Y_-(z) \begin{pmatrix}
                   1 & z^{-3n-2} e^{-3nz} & 0 \\
                   0 & 1 & 0 \\
                   0 & z^{-3n-2}e^{3nz} & 1
                   \end{pmatrix}, \qquad z \in \Gamma,
\end{equation}
where $Y_+(z)$ and $Y_-(z)$ denote the limiting values of $Y(z')$ as $z'$
approaches $z \in \Gamma$ from the inside and outside of $\Gamma$,
respectively.
\item[(3)] For large $z$,
\begin{equation}  \label{eq:Yasym}
    Y(z) = \left( I + \O\left(\frac1z\right) \right)
           \begin{pmatrix}
            z^{n+1} & 0 & 0 \\
            0 & z^{-2n-2} & 0 \\
            0 & 0 & z^{n+1}
            \end{pmatrix}, \qquad z \to \infty.
\end{equation}
\end{enumerate}
The Riemann-Hilbert problem for $Y$ has a unique solution and
$Y_{21}(z) = P_n(z)$, $Y_{22}(z) = z^{-3n-2} Q_n(z)$ for
$z$ outside $\Gamma$, $Y_{22}(z) = z^{-3n-2} E_n(z)$ for
$z$ inside $\Gamma$, and $Y_{23}(z) = R_n(z)$.

A similar characterization holds for the
type II Hermite-Pad\'e approximants. Consider the
following Riemann-Hilbert problem: find a $3\times 3$ matrix valued function
$X: \mathbb{C}\setminus \Gamma \to \mathbb{C}^{3\times 3}$
such that
\begin{enumerate}
\item[(1)] $X$ is analytic in $\mathbb{C} \setminus \Gamma$.
\item[(2)] $X$ satisfies the jump condition
\begin{equation}  \label{eq:Xjump}
   X_+(z) = X_-(z) \begin{pmatrix}
                   1 & 0 & 0 \\
                   -z^{-3n-2} e^{-3nz} & 1 & -z^{-3n-2}e^{3nz} \\
                   0 & 0 & 1
                   \end{pmatrix}, \qquad z \in \Gamma.
\end{equation}
\item[(3)] For large $z$,
\begin{equation}  \label{eq:Xasym}
    X(z) = \left( I + \O\left(\frac1z\right) \right)
           \begin{pmatrix}
            z^{-n-1} & 0 & 0 \\
            0 & z^{2n+2} & 0 \\
            0 & 0 & z^{-n-1}
            \end{pmatrix}, \qquad z \to \infty.
\end{equation}
\end{enumerate}

Comparing the two Riemann-Hilbert problems it is easy to
check that the two solutions $X$ and $Y$ are related by
\begin{equation} \label{relXY}
    X(z) = \left( Y^{-1}(z) \right)^t.
\end{equation}
The matrix $X(z)$ contains the type II polynomials (\ref{typeIIpolys}).
We have the following result.
\begin{theorem}
Let $A_n$, $B_n$, $C_n$, $E_{n}^{(1)}$, and $E_n^{(2)}$ be as
in {\rm (\ref{typeIIpolys})} and {\rm (\ref{error12})}, where we assume that
$A_n$ is normalized to be a monic polynomial of degree $2n+2$.
Then the solution of the above Riemann--Hilbert problem for $X$
is unique and is given by
\begin{equation} \label{eq:Xout}
 X(z) =
         \begin{pmatrix}
         z^{-3n-2} b_{n,n,n+1}(3nz) &  a_{n,n,n+1}(3nz) & z^{-3n-2}
c_{n,n,n+1}(3nz) \\[10pt]
         z^{-3n-2} B_{n}(z) & A_{n}(z) & z^{-3n-2} C_{n}(z) \\[10pt]
         z^{-3n-2} b_{n,n+1,n}(3nz) &  a_{n,n+1,n}(3nz) & z^{-3n-2}
c_{n,n+1,n}(3nz)
         \end{pmatrix},
\end{equation}
for $z$ outside  $\Gamma$, and
\begin{equation} \label{eq:Xin}
 X(z) =  \begin{pmatrix}
         -z^{-3n-2} e_{n,n,n+1}^{(1)}(3nz) &  a_{n,n,n+1}(3nz) & -z^{-3n-2}
e_{n,n,n+1}^{(2)}(3nz) \\[10pt]
         -z^{-3n-2} E_{n}^{(1)}(z) & A_{n}(z) & -z^{-3n-2}
E_{n}^{(2)}(z) \\[10pt]
         -z^{-3n-2} e_{n,n+1,n}^{(1)}(3nz) &  a_{n,n+1,n}(3nz) & -z^{-3n-2}
e_{n,n+1,n}^{(2)}(3nz)
         \end{pmatrix},
\end{equation}
for  $z$  inside $\Gamma$.
In the first rows of {\rm (\ref{eq:Xout})} and {\rm (\ref{eq:Xin})} we use the
Hermite--Pad\'e polynomials of indices
$n,n,n+1$ normalized so that $b_{n,n,n+1}(3nz)$ is a monic polynomial,
and in the third row we use the Hermite--Pad\'e polynomials of indices
$n,n+1,n$ normalized so that $c_{n,n+1,n}(3nz)$ is monic. In addition, we use
\begin{equation} \label{eq:en1n2n3}
    e_{n_1,n_2,n_3}^{(1)}(z) = a_{n_1,n_2,n_3}(z) e^{-z} - b_{n_1,n_2,n_3}(z),
    \quad
    e_{n_1,n_2,n_3}^{(2)}(z) = a_{n_1,n_2,n_3}(z) e^{z} - c_{n_1,n_2,n_3}(z).
    \end{equation}
\end{theorem}
\begin{proof}
The entries in the first and third column of (\ref{eq:Xin}) are analytic at $z=0$
because of the behavior (\ref{typeIIcond}) as $z \to 0$. So $X$ defined
by (\ref{eq:Xout}) and (\ref{eq:Xin}) is analytic on $\mathbb C \setminus \Gamma$.
The jump condition (\ref{eq:Xjump}) is easy to check with the aid of
(\ref{error12}) and (\ref{eq:en1n2n3}). The asymptotic condition (\ref{eq:Xasym})
follows from the degree conditions on the polynomials and the fact that
$b_{n,n,n+1}(3nz)$, $A_n(z)$, and $c_{n,n+1,n}(3nz)$ are monic polynomials.
So $X$ is a solution the Riemann-Hilbert problem.

The uniqueness of the solution is a consequence of Liouville's theorem for
entire functions. Alternatively, it also follows from (\ref{relXY}) and
the fact that the solution of the Riemann-Hilbert problem for $Y$ is
unique, see \cite[Theorem 5.1]{KVW}.
\end{proof}

The fact that the Riemann-Hilbert problems for type I and type II
approximants correspond to each other by applying the inverse transpose
is another way to express the well-known Mahler relations between type
I and type II Hermite-Pad\'e approximants, see \cite{Ma} or \cite[\S 2.1]{NUT}.
This simple relation of taking inverse transpose allows one
to get the asymptotic results for type II approximants directly from the
asymptotic analysis that was made in \cite{KVW} for the solution of
the corresponding type I Riemann-Hilbert problem. In the next section, we
introduce several definitions and notations, some taken from \cite{KVW}
and some new. Then, we state the results.

\section{Asymptotics of type II quadratic Hermite-Pad\'e for the
exponential function}
\subsection{The Riemann surface and related objects}
We define $\mathcal{R}$ as the Riemann surface for the function
\begin{equation}  \label{eq:zw}
    z = z(w)
     = \frac13 \left( \frac1w + \frac{1}{w-1}+\frac{1}{w+1} \right)
     = \frac{w^2-\frac13}{w(w^2-1)}.
\end{equation}
The choice of this particular Riemann surface is motivated by the integral
expressions for the Hermite-Pad\'e approximants to exponentials and
the study of their main contributions by the classical saddle point
method,
see \cite[\S 2.2]{stahl2} or \cite[\S 2.1]{KVW}
for details.
\unitlength=1.00mm
\begin{figure}[t]
\centering
\begin{picture}(140.00,80.00)(0,60)
\put(30.00,140.00){\line(1,0){110.00}}
\put(140.00,140.00){\line(-3,-2){30.00}}
\put(110.00,120.00){\line(-1,0){110.00}}
\put(0.00,120.00){\line(3,2){30.00}}
\put(30.00,110.00){\line(1,0){110.00}}
\put(140.00,110.00){\line(-3,-2){30.00}}
\put(110.00,90.00){\line(-1,0){110.00}}
\put(0.00,90.00){\line(3,2){30.00}}
\put(30.00,80.00){\line(1,0){110.00}}
\put(140.00,80.00){\line(-3,-2){30.00}}
\put(110.00,60.00){\line(-1,0){110.00}}
\put(0.00,60.00){\line(3,2){30.00}}
\bezier{60}(80.00,65.00)(87.00,69.00)(90.00,75.00)
\bezier{60}(50.00,125.00)(54.00,132.00)(60.00,135.00)
\bezier{60}(50.00,95.00)(54.00,102.00)(60.00,105.00)
\bezier{60}(80.00,95.00)(87.00,99.00)(90.00,105.00)
\put(50.00,125.00){\line(0,-1){3.00}}
\put(50.00,120.00){\line(0,-1){3.00}}
\put(50.00,115.00){\line(0,-1){3.00}}
\put(50.00,110.00){\line(0,-1){3.00}}
\put(50.00,105.00){\line(0,-1){3.00}}
\put(50.00,100.00){\line(0,-1){3.00}}
\put(50.00,125.00){\circle*{2.00}}
\put(50.00,95.00){\circle*{2.00}}
\put(60.00,135.00){\line(0,-1){3.00}}
\put(60.00,130.00){\line(0,-1){3.00}}
\put(60.00,125.00){\line(0,-1){3.00}}
\put(60.00,120.00){\line(0,-1){3.00}}
\put(60.00,115.00){\line(0,-1){3.00}}
\put(60.00,110.00){\line(0,-1){3.00}}
\put(60.00,135.00){\circle*{2.00}}
\put(60.00,105.00){\circle*{2.00}}
\put(80.00,95.00){\line(0,-1){3.00}}
\put(80.00,90.00){\line(0,-1){3.00}}
\put(80.00,85.00){\line(0,-1){3.00}}
\put(80.00,80.00){\line(0,-1){3.00}}
\put(80.00,75.00){\line(0,-1){3.00}}
\put(80.00,70.00){\line(0,-1){3.00}}
\put(80.00,95.00){\circle*{2.00}}
\put(80.00,65.00){\circle*{2.00}}
\put(90.00,105.00){\line(0,-1){3.00}}
\put(90.00,100.00){\line(0,-1){3.00}}
\put(90.00,95.00){\line(0,-1){3.00}}
\put(90.00,90.00){\line(0,-1){3.00}}
\put(90.00,85.00){\line(0,-1){3.00}}
\put(90.00,80.00){\line(0,-1){3.00}}
\put(90.00,105.00){\circle*{2.00}}
\put(90.00,75.00){\circle*{2.00}}
\put(56.00,138.00){\makebox(0,0)[cc]{$z_1$}}
\put(46.00,128.00){\makebox(0,0)[cc]{$z_2$}}
\put(84.00,63.00){\makebox(0,0)[cc]{$z_3$}}
\put(94.00,73.00){\makebox(0,0)[cc]{$z_4$}}
\put(84.00,73.00){\makebox(0,0)[cc]{$\Gamma_R$}}
\put(55.00,127.00){\makebox(0,0)[cc]{$\Gamma_P$}}
\put(127.00,126.00){\makebox(0,0)[cc]{$\mathcal{R}_P$}}
\put(127.00,96.00){\makebox(0,0)[cc]{$\mathcal{R}_Q$}}
\put(127.00,66.00){\makebox(0,0)[cc]{$\mathcal{R}_R$}}
\end{picture}
\caption{The Riemann surface $\mathcal{R}$}
\label{fig:surface}
\end{figure}
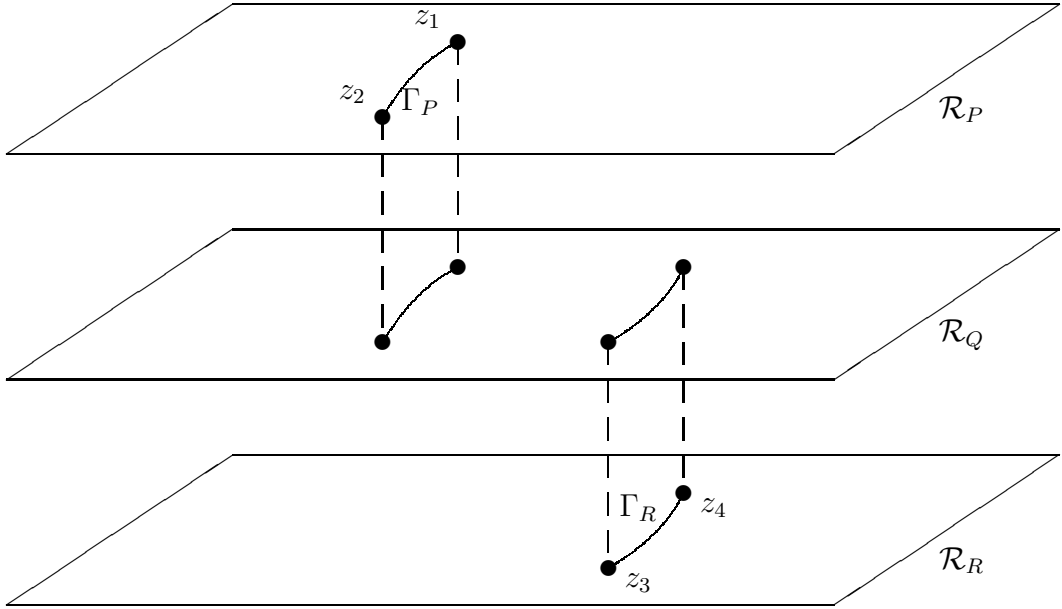

The rational function (\ref{eq:zw}) has three inverse mappings. These are
the three solutions of the cubic equation
\begin{equation}  \label{eq:cubic}
   zw^3-w^2-zw+\frac{1}{3} = 0.
\end{equation}
The Riemann surface $\mathcal{R}$ consists of three sheets $\mathcal{R}_P$,
$\mathcal{R}_Q$, and $\mathcal{R}_R$ as shown in Figure \ref{fig:surface}. The
bijective mapping
$\psi : \mathcal{R} \to \overline{\mathbb{C}}$ is the inverse of (\ref{eq:zw}).
We denote its restriction to the three sheets by $\psi_P$, $\psi_Q$,
and $\psi_R$, respectively. So $\psi_P(z)$, $\psi_Q(z)$, and $\psi_R(z)$
are the three solutions of (\ref{eq:cubic}), and we assume that the
choice of the three sheets is such that
$\psi_P(\infty) = -1$, $\psi_Q(\infty) = 0$, and $\psi_R(\infty) = 1$.
The Riemann surface has four branch points
$z_1=z(w_1)$, $z_2=z(w_2)$, $z_3=z(w_3)$, $z_4=z(w_4)$, which
are related to the points $w_1$, $w_2$, $w_3$, $w_4$ for which $z'(w) = 0$.
Simple calculations give
\begin{equation}  \label{eq:wk}
  w_k = 3^{-1/4} \omega_8^{-2k-1}, \qquad k=1,2,3,4,
\end{equation}
where $\omega_8 = e^{2\pi i/8}$  is the primitive 8th root of
unity. The corresponding values of $z_k = z(w_k)$ are
\begin{equation} \label{eq:zk}
z_1 = 3^{-1/4} \omega_{24}^{7},\quad
z_2 = 3^{-1/4} \omega_{24}^{17},\quad
z_3 = 3^{-1/4} \omega_{24}^{19},\quad
z_4 = 3^{-1/4} \omega_{24}^{5},
\end{equation}
where $\omega_{24} = e^{2\pi i/24}$ is the primitive 24th root of unity.

\begin{figure}[tb]
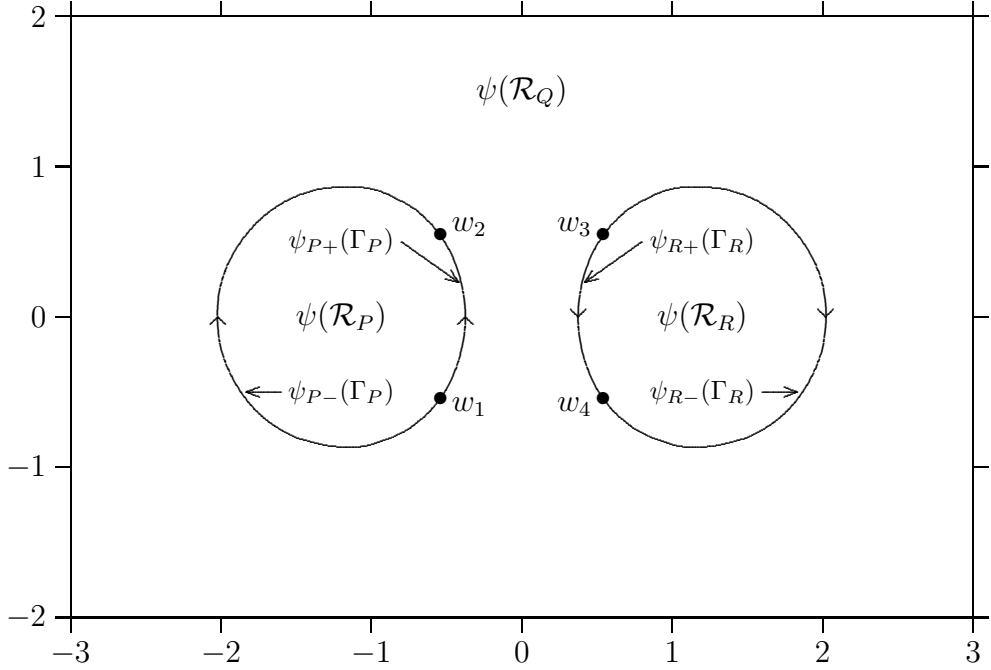

\hfil
\input traject8
\hfil
\caption{$\psi$-image of the Riemann surface $\mathcal{R}$}  \label{fig:psi}
\end{figure}

The sheets $\mathcal{R}_P$ and $\mathcal{R}_Q$ are glued together
along a cut $\Gamma_P$
connecting two branch points $z_1$ and $z_2$ in the left half-plane,
and the sheets
$\mathcal{R}_Q$ and
$\mathcal{R}_R$ are glued together along a cut $\Gamma_R$ connecting the
other two branch points $z_3$ and $z_4$, in the right
half-plane. Moreover, these cuts are chosen so as to satisfy
\begin{equation} \label{eq:GammaP}
    \frac{3}{2\pi i} \int_{z_1}^z (\psi_Q-\psi_P)_+(s) ds \in \mathbb R,\quad
z\in \Gamma_P,
\end{equation}
\begin{equation}  \label{eq:GammaR}
  \frac{3}{2\pi i} \int_{z_3}^z (\psi_Q-\psi_R)_+(s)\,ds \in \mathbb{R},\quad z\in \Gamma_{R}.
\end{equation}
In particular, $\Gamma_P$ and $\Gamma_R$ are trajectories of the quadratic
differentials $-(\psi_Q-\psi_P)^2(z)dz^2$ and $-(\psi_Q-\psi_R)^2(z)dz^2$
respectively, see \cite{POM, STR}. The $\psi$-images of these curves are
shown in Figure \ref{fig:psi}.

\begin{figure}[tb]
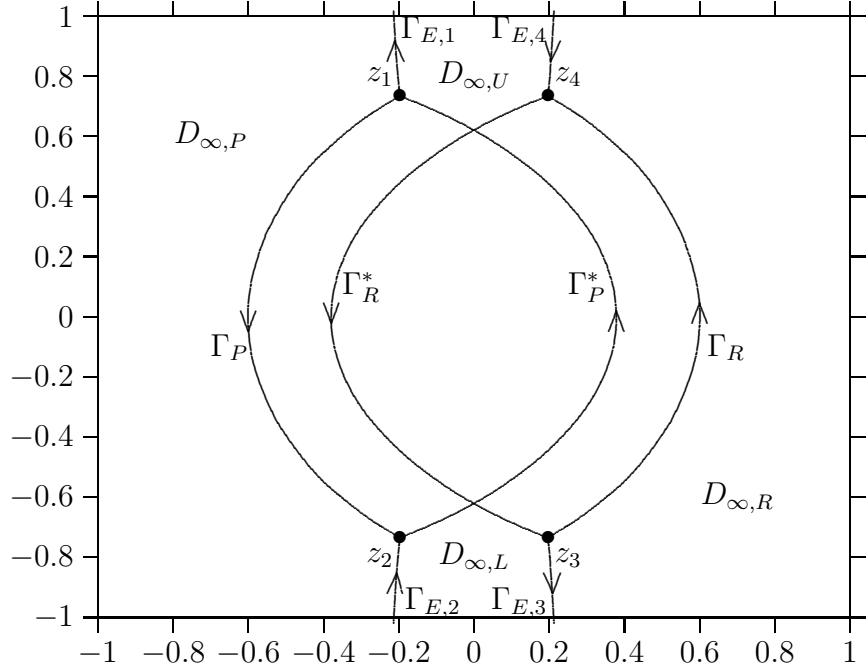

\hfil
\input traject
\hfil
\caption{Curves for which either $\frac{3}{2\pi i} \int_{z_1}^z (\psi_Q-\psi_P)(s)\,ds$
or $\frac{3}{2\pi i} \int_{z_3}^z (\psi_Q-\psi_R)(s)\,ds$ is real.}
 \label{fig:traject}
\end{figure}

Actually, there are four analytic curves such that (\ref{eq:GammaP}),
resp. (\ref{eq:GammaR}), holds. We denote them by
$\Gamma_{P}$, $\Gamma_{P}^{*}$, $\Gamma_{E,1}$, $\Gamma_{E,2}$,
respectively
$\Gamma_{R}$, $\Gamma_{R}^{*}$, $\Gamma_{E,3}$, $\Gamma_{E,4}$, see
Figure \ref{fig:traject}.
The particular geometry of these curves may be derived from the local and
global properties of trajectories of quadratic differentials, see
\cite[Proposition 2.2]{KVW}. These curves divide the complex plane
into different domains. We denote the unbounded domains by
$D_{\infty,P}$, $D_{\infty,R}$,
$D_{\infty,U}$ and $D_{\infty,L}$, as shown in Figure
\ref{fig:traject}. We denote by $D_P^*$ the bounded domain
bounded by $\Gamma_P$ and $\Gamma_P^*$, and by $D_R^*$, its mirror
image with respect to the imaginary axis.

\subsection{Definition of measures and weak convergence of zeros}
We now define a measure on each of the previous curves. The complex
line element $ds$ is defined according to
the orientation of these curves given in Figure \ref{fig:traject}.
\begin{definition}
We define a measure $\mu_P$ on $\Gamma_P$ by
\begin{equation} \label{eq:muP}
    d\mu_P(s) = \frac{3}{2\pi i} (\psi_Q-\psi_P)_+(s)\, ds \qquad s
    \in \Gamma_P,
\end{equation}
and a measure $\mu_R$ on $\Gamma_R$ by
\begin{equation} \label{eq:muR}
  d\mu_R(s) = \frac{3}{2\pi i} (\psi_Q-\psi_R)_+(s)\, ds \qquad s \in \Gamma_R.
\end{equation}
The measure $\mu_A$ on $\Gamma_P\cup\Gamma_R$ is defined by
\begin{equation} \label{eq:muA}
    \mu_A = \mu_P + \mu_R.
\end{equation}
We define a measure $\mu_B$ on $\Gamma_P^{*}$ by
\begin{equation} \label{eq:muB}
    d\mu_{B}(s) =
        \frac{3}{2\pi i} (\psi_Q-\psi_P)(s)\, ds \qquad s \in \Gamma_P^*,
\end{equation}
a measure $\mu_{C}$ on $\Gamma_R^{*}$ by
\begin{equation} \label{eq:muC}
    d\mu_{C}(s) =
        \frac{3}{2\pi i} (\psi_Q-\psi_R)(s)\, ds \qquad s \in \Gamma_R^*,
\end{equation}
The measures $\mu_{E^{(1)}}$ on $\Gamma_{E,1} \cup \Gamma_{E,2}$ and
$\mu_{E^{(2)}}$ on $\Gamma_{E,3} \cup \Gamma_{E,4}$ are defined as
\begin{eqnarray}  \label{eq:muE1}
  d\mu_{E^{(1)}}(s) & = &
       \frac{3}{2\pi i} (\psi_Q-\psi_P)(s)\, ds \qquad s \in
\Gamma_{E,1} \cup \Gamma_{E,2}, \\
      \label{eq:muE2}
  d\mu_{E^{(2)}}(s) & = &
       \frac{3}{2\pi i} (\psi_Q-\psi_R)(s)\, ds \qquad s \in
\Gamma_{E,3} \cup \Gamma_{E,4}.
         \end{eqnarray}
\end{definition}

The facts that $\mu_P$ and $\mu_R$ are probability measures
and that $\mu_{E^{(1)}}$ and $\mu_{E^{(2)}}$ are positive measures of
infinite masses are proven in \cite[Theorem 2.4]{KVW}. So $\mu_A$
is a positive measure with total mass $2$. A similar proof shows
that $\mu_B$ and $\mu_C$ are positive measures of total mass $2$,
see also the first part of the proof of \cite[Proposition 4.8]{KVW}.

Denote by $\nu_{A_{n}}$, $\nu_{B_{n}}$, $\nu_{C_{n}}$ the normalized zero
counting measures of the polynomials $A_{n}$, $B_{n}$, $C_{n}$. Furthermore,
define the zero
counting measures for the remainder functions $E_n^{(1)}$ and
$E_n^{(2)}$ as
\[ \nu_{E_n^{(1)}} = \frac{1}{2n} \sum_{\stackrel{E_n^{(1)}(z) =
0}{z\neq 0}} \delta_{z}, \qquad
\nu_{E_n^{(2)}} = \frac{1}{2n} \sum_{\stackrel{E_n^{(2)}(z) =
0}{z\neq 0}} \delta_{z},\]
where
the normalization by $2n$ now corresponds to the degree of approximation
and the $3n+2$ interpolatory zeros of $E_n^{(1)}$ and $E_n^{(2)}$ at
the origin have been
excluded. With all of the previous measures, we can state our
first result.
\begin{theorem}
\label{lim-zeros}
We have
\begin{equation} \label{eq:convnuPQRn}
    \nu_{A_n} \stackrel{*}{\to} \frac{1}{2} \mu_A, \qquad
   \nu_{B_n} \stackrel{*}{\to} \frac{1}{2} \mu_B, \qquad
   \nu_{C_n} \stackrel{*}{\to} \frac{1}{2} \mu_C,
\end{equation}
where the convergence is in the sense of weak$^*$ convergence of measures,
i.e., $\mu_n \stackrel{*}{\to} \mu$ if $\int f d\mu_n \to \int fd\mu$ for
every bounded continuous function $f$.
Furthermore, we have
\begin{equation} \label{eq:convnuEn}
    \nu_{E_n^{(1)}} \stackrel{*}{\to} \mu_{E^{(1)}},\qquad
\nu_{E_n^{(2)}} \stackrel{*}{\to} \mu_{E^{(2)}},
\end{equation}
in the sense that
\[ \lim_{n\to\infty} \int f(s) d\nu_{E_n^{(i)}}(s) = \int f(s)
d\mu_{E^{(i)}}(s),\quad i=1,2, \]
for every continuous function $f$ such that $f(s) = \O(s^{-2})$
as $s \to \infty$.
\end{theorem}

\subsection{Definition of $g$-functions}
Next, we introduce the complex logarithmic potentials of the
measures $\mu_P$, $\mu_Q$, $\mu_{A}$, $\mu_{B}$, and $\mu_{C}$.
\begin{definition}
We set for $S \in \{P,Q,A,B,C\}$,
\begin{eqnarray}
  g_S(z) & = & \int \log (z-s) \, d\mu_S(s),   \qquad z \in \mathbb C \setminus \supp(\mu_S), \label{eq:gS}
\end{eqnarray}
which is defined modulo $2\pi i$.
\end{definition}
We will also need two more functions, namely,
\begin{equation}  \label{eq:varphiP}
   \varphi_P(z) = \frac{3}{2} \int_{z_1}^z (\psi_Q-\psi_P)(s)\, ds,
\end{equation}
and
\begin{equation} \label{eq:varphiR}
   \varphi_R(z) = \frac{3}{2} \int_{z_3}^z (\psi_Q-\psi_R)(s)\, ds,
\end{equation}
where the paths of integration are in $\C
\setminus\left(\Gamma_{P}\cup\Gamma_{R}\cup\{0 \} \right)$.
These functions were used in \cite{KVW} and we refer to \cite[Lemma
2.7]{KVW} for their main properties. Note that the logarithmic
potentials in (\ref{eq:gS}) have jumps on the support
of their measures which can be expressed in terms of $\varphi_{P}$ and
$\varphi_{R}$, namely we have
\begin{equation}\label{jumpAP}
g_{A+} (z)-g_{A-} (z)=-2\varphi_{P+} (z)=2\varphi_{P-} (z),\quad z\in
\Gamma_{P},
\end{equation}
\begin{equation}\label{jumpAR}
g_{A+} (z)-g_{A-} (z)=-2\varphi_{R+} (z)=2\varphi_{R-} (z),\quad z\in
\Gamma_{R},
\end{equation}
\begin{equation}\label{jumpB}
g_{B+} (z)-g_{B-} (z)=-2\varphi_{P} (z),\quad z\in
\Gamma_{P}^{*},
\end{equation}
\begin{equation}\label{jumpC}
g_{C+} (z)-g_{C-} (z)=-2\varphi_{R} (z),\quad z\in
\Gamma_{R}^{*}.
\end{equation}
\subsection{Strong asymptotics away from the zeros}

As in \cite{KVW}, we will use the function
$\sqrt{3w^4+1}$ which branches at the four points $w_k$ given in (\ref{eq:wk}).
We choose as cuts for
this function the two curves $\psi_{P+}(\Gamma_P)$ and
$\psi_{R+}(\Gamma_R)$ (see Figure \ref{fig:psi}), and assume that it
is positive for large positive $w$. So, in particular we have
that $\sqrt{3w^4+1} = -1$ for $w=0$.

The following theorem gives the strong asymptotics of the polynomials
$A_n$, $B_n$, $C_n$ and the remainder terms $E_n^{(1)}$ and
$E_{n}^{(2)}$ away from their zeros.
\begin{theorem}
\label{asymp}
With the functions defined above, we have
\begin{equation} \label{asympA}
A_n(z) = -\frac{z^2 (\psi_{Q}^{2} (z) -1)^{2} e^{ng_A(z)}}
    {\sqrt{3\psi_Q^4(z)+1}}
    \left(1+ {\O}\left(\frac{1}{n}\right)\right)
\end{equation}
uniformly for $z$ in compact subsets of $\mathbb C \setminus
(\Gamma_P\cup \Gamma_{R})$,
\begin{equation} \label{asympB}
B_n(z) = \left\{ \begin{array}{ll} \ds
    -(-2)^{n}\frac{z^2 (\psi_{P}^{2} (z) -1)^{2}e^{n g_B(z)}}
    {\sqrt{3\psi_P^4(z) + 1}}
\left(1 + {\O}\left(\frac{1}{n}\right)\right)
    & \mbox{ for } z \in \mathbb C \setminus
(D_P^{*}\cup \Gamma_{P}^{*}), \\[10pt]
    \ds
    -(-2)^{n}\frac{z^2 (\psi_{Q}^{2} (z) -1)^{2} e^{n g_B(z)}}
    {\sqrt{3\psi_Q^4(z) + 1}}
\left(1 + {\O}\left(\frac{1}{n}\right)\right)
    & \mbox{ for } z \in D_{P}^{*},
    \end{array} \right.
\end{equation}
uniformly for $z$ in compact subsets of $\mathbb C \setminus
\Gamma_P^{*}$, and
\begin{equation} \label{asympC}
C_n(z) = \left\{ \begin{array}{ll} \ds
    -(-2)^{n}\frac{z^2 (\psi_{R}^{2} (z) -1)^{2}e^{n g_C(z)}}
    {\sqrt{3\psi_R^4(z) + 1}} \left(1 + {\O}\left(\frac{1}{n}\right)\right)
    & \mbox{ for } z \in \mathbb C \setminus
(D_R^{*}\cup \Gamma_{R}^{*}), \\[10pt]
    \ds
    -(-2)^{n}\frac{z^2 (\psi_{Q}^{2} (z) -1)^{2}e^{n g_C(z)}}
    {\sqrt{3\psi_Q^4(z) + 1}}
    \left(1 + {\O}\left(\frac{1}{n}\right)\right)
        & \mbox{ for } z \in D_{R}^{*}, \\[10pt]
    \end{array} \right.
\end{equation}
uniformly for $z$ in compact subsets of $\mathbb C \setminus
\Gamma_R^{*}$.
Furthermore we have
\begin{equation} \label{asympE1}
 E_n^{(1)}(z) = \left\{ \begin{array}{lll} \ds
    - \frac{z^2 (\psi_{Q}^{2} (z) -1)^{2}e^{n(g_{A} (z) -3z)}}
    {\sqrt{3\psi_Q^4(z)+1}}
    \left(1 + {\O}\left(\frac{1}{n}\right)\right)
    & \mbox{ for } z \in D_{\infty,P}, \\[10pt] \ds
    (-2)^{n}\frac{z^2(\psi_{P}^{2} (z) -1)^{2}e^{ng_B(z)}}
    {\sqrt{3\psi_P^4(z) + 1}}  \left(1 + {\O}\left(\frac{1}{n}\right)\right)
    & \mbox{ for } z \in \C \setminus (D_{\infty,P} \cup D_{P}^{*}\\[-10pt]
    & \qquad \quad \cup
\Gamma_{P}\cup
\Gamma_{E,1}\cup\Gamma_{E,2}),\\[10pt] \ds
    (-2)^{n}\frac{z^{3n+2}(\psi_{P}^{2} (z) -1)^{2} e^{-n g_P(z)}}
    {\sqrt{3\psi_P^4(z) + 1}}  \left(1 + {\O}\left(\frac{1}{n}\right)\right)
    & \mbox{ for } z \in D_{P}^{*},
    \end{array} \right.
\end{equation}
uniformly for $z$ in compact subsets of $\mathbb{C} \setminus (
\Gamma_{E,1}\cup\Gamma_{E,2})$, and
\begin{equation} \label{asympE2}
 E_n^{(2)}(z) = \left\{ \begin{array}{ll} \ds
    - \frac{z^2(\psi_{Q}^{2} (z) -1)^{2}e^{n (g_{A} (z) +3z)}}
    {\sqrt{3\psi_Q^4(z)+1}}
    \left(1 + {\O}\left(\frac{1}{n}\right)\right)
    & \mbox{ for } z \in D_{\infty,R}, \\[10pt] \ds
    (-2)^{n}\frac{z^2(\psi_{R}^{2} (z) -1)^{2}e^{ng_C(z)}}
    {\sqrt{3\psi_R^4(z) + 1}}  \left(1 + {\O}\left(\frac{1}{n}\right)\right)
    & \mbox{ for } z \in \C \setminus (D_{\infty,R} \cup D_{R}^{*}\\[-10pt]
    & \qquad \quad \cup \Gamma_{R}\cup
    \Gamma_{E,3}\cup\Gamma_{E,4}),\\[10pt] \ds
    (-2)^{n}\frac{z^{3n+2}(\psi_{R}^{2} (z) -1)^{2}e^{-n g_R(z)}}
    {\sqrt{3\psi_R^4(z) + 1}}  \left(1 + {\O}\left(\frac{1}{n}\right)\right)
    & \mbox{ for } z \in D_{R}^{*},
    \end{array} \right.
\end{equation}
uniformly for $z$ in compact subsets of $\mathbb{C} \setminus (
\Gamma_{E,3}\cup\Gamma_{E,4})$.
\end{theorem}

\subsection{Asymptotics near the curves $\Gamma_{P}$,
$\Gamma_{P^{*}}$, $\Gamma_{R}$, $\Gamma_{R^{*}}$, and
$\Gamma_{E^{(i)}}$, $i=1,2,3,4$.}
\begin{figure}[tb]
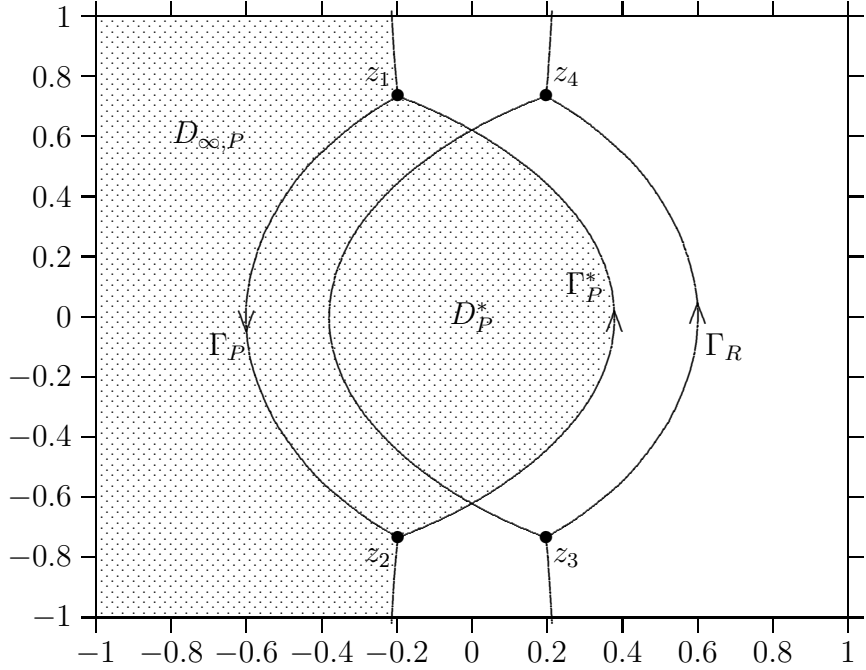

\hfil
\input traject1
\hfil \caption{The shaded region is where $\Re \varphi_P$ is
negative. It consists of the two parts $D_{\infty, P}$ and
$D_P^*$, where $D_P^*$ is bounded by $\Gamma_P$ and $\Gamma_P^*$.}
\label{fig:Re}
\end{figure}

We state in the next theorem uniform asymptotics which are valid in
neighborhoods of the curves where the zeros of the polynomials $A_{n}$, $B_{n}$,
$C_{n}$ and the error functions $E_{n}^{(1)}$ and $E_{n}^{(2)}$
accumulate, see Theorem \ref{lim-zeros}.
\begin{theorem}
\label{asym-curves}
Uniformly for $z$ in compact subsets of the region
$\left(\Gamma_P \setminus \{z_1, z_2 \} \right) \cup
    D_{\infty,P} \cup D_P^*$
(this is the shaded region in Figure {\rm \ref{fig:Re}}),
we have
\begin{equation} \label{asympcurA1}
A_n(z)= e^{n g_A(z)}
    \left[ -\frac{z^2(\psi_{Q}^{2}(z) -1)^{2}}{\sqrt{3\psi_{Q}^4(z)+1}} \left(1 +
    \O\left(\frac{1}{n}\right)\right)
        \mp \frac{z^2(\psi_{P}^{2} (z) -1)^{2}e^{2n \varphi_P(z)}}
        {\sqrt{3\psi_{P}^4(z)+1}}
    \left(1 + \O\left(\frac{1}{n}\right)\right) \right]
\end{equation}
where the $+$sign holds in $D_{\infty,P}$ and the $-$sign holds in $D_{P}^{*}$.

Uniformly for $z$ in compact subsets of the region
$\left(\Gamma_R \setminus \{z_3, z_4 \} \right) \cup
    D_{\infty,R} \cup D_R^*$
we have
\begin{equation} \label{asympcurA2}
A_n(z)= e^{n g_A(z)}
    \left[- \frac{z^2(\psi_{Q}^{2} (z) -1)^{2}}
    {\sqrt{3\psi_{Q}^4(z)+1}} \left(1 +
        \O\left(\frac{1}{n}\right)\right)
        \mp \frac{z^2(\psi_{R}^{2} (z) -1)^{2}e^{2n \varphi_R(z)}}
        {\sqrt{3\psi_{R}^4(z)+1}}
    \left(1 + \O\left(\frac{1}{n}\right)\right) \right]
\end{equation}
where the $+$sign holds in $D_{\infty,R}$ and the $-$sign holds in $D_{R}^{*}$.

Uniformly for $z$ in compact subsets of
$\C \setminus \overline{D_{\infty,P}}$,
we have
\begin{multline} \label{asympcurB}
B_n(z)= - (-2)^{n} e^{n (g_B(z)\pm\varphi_{P} (z))}
    \left[ \frac{z^2(\psi_{Q}^{2} (z) -1)^{2}e^{-n \varphi_P(z)}}
    {\sqrt{3\psi_{Q}^4(z)+1}} \left(1 +
    \O\left(\frac{1}{n}\right)\right) \right. \\   \left.
        + \frac{z^2(\psi_{P}^{2} (z) -1)^{2}e^{n \varphi_P(z)}}
        {\sqrt{3\psi_{P}^4(z)+1}}
        \left(1 + \O\left(\frac{1}{n}\right)\right) \right]
\end{multline}
where the $+$sign holds in $D_{P}^{*}$ and the $-$sign holds in
$\C \setminus\left(\overline{D_{\infty,P}\cup D_{P}^{*}}\right)$.

Uniformly for $z$ in compact subsets of
$\C \setminus \overline{D_{\infty,R}}$,
we have
\begin{multline} \label{asympcurC}
C_n(z)= - (-2)^{n} e^{n (g_C(z)\pm\varphi_{R} (z))}
    \left[ \frac{z^2(\psi_{Q}^{2} (z) -1)^{2}e^{-n \varphi_R(z)}}
    {\sqrt{3\psi_{Q}^4(z)+1}} \left(1 +
    \O\left(\frac{1}{n}\right)\right) \right. \\   \left.
        + \frac{z^2(\psi_{R}^{2} (z) -1)^{2}e^{n \varphi_R(z)}}
        {\sqrt{3\psi_{R}^4(z)+1}}
        \left(1 + \O\left(\frac{1}{n}\right)\right) \right]
\end{multline}
where the $+$sign holds in $D_{R}^{*}$ and the $-$sign holds in
$\C \setminus\overline{D_{\infty,R}\cup D_{R}^{*}}$.

Uniformly for $z$ in compact subsets of
$\C \setminus \overline{D_{P}^{*}}$,
we have
\begin{multline} \label{asympcurE1}
E_n^{(1)}(z)= e^{n (g_A(z)-3z)}
    \left[ \frac{z^2(\psi_{P}^{2} (z) -1)^{2}e^{2n \varphi_P(z)}}
    {\sqrt{3\psi_{P}^4(z)+1}} \left(1 +
    \O\left(\frac{1}{n}\right)\right) \right. \\   \left.
        - \frac{z^2(\psi_{Q}^{2} (z) -1)^{2}}
        {\sqrt{3\psi_{Q}^4(z)+1}}
        \left(1 + \O\left(\frac{1}{n}\right)\right) \right].
\end{multline}

Uniformly for $z$ in compact subsets of
$\C \setminus \overline{D_{R}^{*}}$,
we have
\begin{multline} \label{asympcurE2}
E_n^{(2)}(z)= e^{n (g_A(z)+3z)}
    \left[ \frac{z^2(\psi_{R}^{2} (z) -1)^{2}e^{2n \varphi_R(z)}}
    {\sqrt{3\psi_{R}^4(z)+1}} \left(1 +
    \O\left(\frac{1}{n}\right)\right) \right. \\   \left.
        - \frac{(z^2 3\psi_{Q}^{2} (z) -1)^{2}}
        {\sqrt{3\psi_{Q}^4(z)+1}}
        \left(1 + \O\left(\frac{1}{n}\right)\right) \right].
\end{multline}
\end{theorem}

The formulas of Theorem \ref{asym-curves} are valid on the
curves where the zeros accumulate. For example, (\ref{asympcurA1})
holds on $\Gamma_{P}$ away
from the branch points $z_{1}$ and $z_{2}$ and (\ref{asympcurA2})
holds on $\Gamma_{R}$ away
from the branch points $z_{3}$ and $z_{4}$, that is on the two
curves where the zeros of $A_{n}$ accumulate.

It has been proved in \cite[Lemma 2.7]{KVW} that the real part of
$\varphi_{P}$ is negative in $D_{\infty,P}\cup D_{P}^{*}$ and that it
is positive in the remaining part of the complex plane, see Figure
\ref{fig:Re}.
For $z$ away from $\Gamma_P$ so that $\Re \varphi_P(z) < 0$, the
asymptotic formula
(\ref{asympcurA1}) reduces to (\ref{asympA}). On $\Gamma_P$ we have
$\Re \varphi_P(z) =0$,
and then the two terms in (\ref{asympcurA1}) are of comparable
magnitude. For $z \in \Gamma_P$,
(\ref{asympcurA1}) can be re-written as
\[
    A_n(z)= -
    \frac{z^2(\psi_{Q-}^{2} (z) -1)^{2}e^{n g_{A-}(z)}}
    {\sqrt{3\psi_{Q-}^4(z)+1}} \left(1 +
    \O\left(\frac{1}{n}\right)\right)
        -\frac{z^2(\psi_{Q+}^{2} (z) -1)^{2}e^{n g_{A+}(z)}}
        {\sqrt{3\psi_{Q+}^4(z)+1}}
    \left(1 + \O\left(\frac{1}{n}\right)\right).
\]
Similar remarks hold for the other formulas.

\subsection{Asymptotics near the branch points}
As in \cite{KVW} the asymptotic formulas near the branch points involve the
classical Airy function $\Ai$.
which is the solution of the differential equation
$ y''(z)=z y(z)$ that satisfies
$$\Ai(z) = \frac{1}{2\sqrt{\pi}} z^{-1/4} e^{-\frac{2}{3}z^{3/2}}
    \left(1 + \O\left(\frac{1}{z^{3/2}}\right)\right)
$$
as $z \to \infty$ with $| \arg z | < \pi$.
Because of symmetry, we may restrict our study to the asymptotic behavior of
$A_n$, $B_n$ and $E_n^{(1)}$ near $z_1$. Similar results hold true
near the other branch
points.

The asymptotics to be given involve a function $f_1$ that is
defined in \cite[Section 2.7]{KVW} as
\begin{equation} \label{deffz}
    f_1(z) = \left[\frac{3}{2} \varphi_P(z) \right]^{2/3},
\end{equation}
where we take the branch of  $\varphi_P(z)$ which behaves like
\[ \varphi_P(z) = c(z-z_1)^{3/2} + \O\left((z-z_1)^{5/2}\right) \]
as $z \to z_1$, with $c \neq 0$.
The $2/3$rd power in (\ref{deffz}) is taken so that
$f_1(z)$ is real and negative for $z \in \Gamma_P$.
Then explicit calculations show that
\begin{equation} \label{def-f}
    f_1(z) = c_1 (z-z_1) + \O\left((z-z_1)^2\right)
\end{equation}
as $z \to z_1$ with \begin{equation}  \label{c_1}
    c_1 = f_1'(z_1) = 2^{1/3} 3^{5/12} e^{-\frac{7}{36} \pi i}.
\end{equation}
\begin{theorem}
\label{asym-branch}
There is a $\delta > 0$, such that we have, uniformly for $|z-z_1| < \delta$,
\begin{multline} \label{asympbrA}
A_n(z)  =  -iz^{-1}e^{-3z}\sqrt{\pi} e^{(n+1) (g_A(z) + \varphi_P(z))} \\
 \times \left[ n^{1/6} \widetilde{h}_1(z) \Ai\left( (n+1)^{2/3} f_1(z)\right) \left(1 + \O\left(\frac{1}{n}\right)\right)
  \right. \\   \left.
    + n^{-1/6} \widetilde{h}_2(z) \Ai'\left((n+1)^{2/3} f_1(z)\right) \left(1+\O\left(\frac{1}{n} \right)\right)
\right],
\end{multline}
\begin{multline} \label{asympbrB}
B_n(z) =  z^{-1}e^{-i\pi/6}\sqrt{\pi} e^{(n+1) (g_A(z) + \varphi_P(z)-3z)} \\
    \times \left[ n^{1/6} \widetilde{h}_1(z)
    \Ai\left(e^{-2\pi i/3} (n+1)^{2/3} f_1(z)\right)
    \left(1 + \O\left(\frac{1}{n}\right)\right)
    \right. \\ \left.
    + n^{-1/6} \widetilde{h}_2(z) e^{-2\pi i/3}
\Ai'\left(e^{-2\pi i/3} (n+1)^{2/3} f_1(z)\right)
\left(1+\O\left(\frac{1}{n} \right)\right)
\right],
\end{multline}
and
\begin{multline} \label{asympbrE1}
E_n^{(1)}(z) = z^{-1}e^{-i\pi/6}\sqrt{\pi} e^{(n+1) (g_A(z) + \varphi_P(z)-3z)} \\
    \times \left[ n^{1/6} \widetilde{h}_1(z) e^{-2 \pi i/3}
\Ai\left(e^{2\pi i/3} (n+1)^{2/3} f_1(z)\right)
    \left(1 + \O\left(\frac{1}{n}\right)\right)
    \right. \\ \left.
    + n^{-1/6} \widetilde{h}_2(z)
\Ai'\left(e^{2\pi i/3} (n+1)^{2/3} f_1(z)\right)
    \left(1+\O\left(\frac{1}{n} \right)\right)
\right],
\end{multline}
where $\widetilde{h}_1$ and $\widetilde{h}_2$ are two analytic functions without zeros
in $|z-z_1| < \delta$, which have explicit expressions
$$    \widetilde{h}_1(z) = \left(\widetilde{N}_{21}(z) + i ze^{3z}
\widetilde{N}_{22}(z) \right)
f_1(z)^{1/4},
$$
with the branch of the fourth root in $f_1(z)^{1/4}$ taken with a cut along $\Gamma_P$, and
$$    \widetilde{h}_2(z) = \left(\widetilde{N}_{21}(z) - ize^{3z}
\widetilde{N}_{22}(z) \right)
f_1(z)^{-1/4}.
$$
Here
$$
\widetilde{N}_{21}(z) = \frac{(\psi_P^2(z)-1)^2e^{g_P(z)}}
{2\sqrt{3 \psi_P^4(z) + 1}}, \qquad
\widetilde{N}_{22}(z) =
- \frac{z^2 (\psi_Q^2(z)-1)^2e^{-g_A(z)}}
{\sqrt{3 \psi_Q^4(z) + 1}}.
$$
\end{theorem}
In Theorem 2.11 of \cite{KVW} we used functions $h_1$, $h_2$, $N_{21}$ and
$N_{22}$
to state asymptotics near the branch point $z_{1}$ for type I Hermite-Pad\'e
approximants. The functions $\widetilde{h}_1$, $\widetilde{h}_2$,
$\widetilde{N}_{21}$ and
$\widetilde{N}_{22}$ are the corresponding functions
for type II Hermite-Pad\'e approximants.
Note that, by (\ref{jumpAP}), the function $g_A + \varphi_P$ is
analytic near $z=z_1$.

From the asymptotics near the branch points, one can deduce the
behavior of the extreme zeros of $A_n$, $B_n$, $C_n$, $E_{n}^{(1)}$
and $E_n^{(2)}$ near the
branch points. We only
state the result for the zeros of $A_n$, $B_n$ and $E_n^{(1)}$ near $z_1$.
Recall that the Airy function $\Ai$ has only negative real zeros,
which we denote by $0 > -\iota_1 > -\iota_2 > \cdots > -\iota_{\nu} > \cdots$.
\begin{corollary} \label{asym-zero}
Let $z^A_{\nu,n}$, $\nu=1,\ldots,2n+2$, be the zeros of $A_n$, ordered
by increasing distance to $z_1$. Then for every $\nu \in \mathbb N$,
we have
\begin{eqnarray} \nonumber
    z^A_{\nu,n} & = & z_1 -  \frac{\iota_\nu}{f_1'(z_1)} n^{-2/3}  +
\O\left(\frac{1}{n}\right)
\\
    & = & \label{zerosPn}
     z_1  +  2^{-1/3} 3^{-5/12} e^{-\frac{29}{36}\pi i}
    \iota_{\nu} n^{-2/3} + \O\left(\frac{1}{n}\right),
\end{eqnarray}
as $n \to \infty$.
Let $z^B_{\nu,n}$, $\nu =1,\ldots, 2n$, be the zeros of $B_n$, ordered
by increasing distance to $z_1$. Then for every $\nu \in \mathbb N$,
\begin{eqnarray} \nonumber
    z^B_{\nu,n} & = & z_1 - e^{2\pi i/3} \frac{\iota_\nu}{f'_1(z_1)}
    n^{-2/3} + \O\left(\frac{1}{n}\right)  \\
    & = & \label{zerosQn}
    z_1 + 2^{-1/3} 3^{-5/12} e^{-\frac{5}{36} \pi i}
    \iota_{\nu} n^{-2/3} + \O\left(\frac{1}{n}\right),
\end{eqnarray}
as $n \to \infty$.
Let $z^{E^{(1)}}_{\nu,n}$, $\nu \geq 1$, be the zeros of $E^{(1)}_n$, ordered
by increasing distance to $z_1$. Then for every $\nu \in \mathbb N$,
\begin{eqnarray} \nonumber
    z^{E^{(1)}}_{\nu,n} & = & z_1 - e^{-2\pi i/3} \frac{\iota_\nu}{f'_1(z_1)}
    n^{-2/3} + \O\left(\frac{1}{n}\right)  \\
    & = & \label{zerosEn}
    z_1 + 2^{-1/3} 3^{-5/12} e^{\frac{19}{36} \pi i}
    \iota_{\nu} n^{-2/3} + \O\left(\frac{1}{n}\right),
\end{eqnarray}
as $n \to \infty$.
\end{corollary}
We note that the previous formulas are exactly the same as those
obtained for the zeros of the Hermite-Pad\'e approximants of type I, $P_{n}$,
$Q_{n}$, and $E_{n}$ near $z_{1}$ in \cite[Corollary 2.12]{KVW}.

\section{Proofs of the asymptotic formulas}
\subsection{Transformations of the Riemann-Hilbert problem}
In \cite{KVW}, the analysis of the Riemann-Hilbert problem for the
matrix $Y$ was made through a sequence of explicit transformations
\begin{equation} \label{transformationsY}
Y \mapsto U \mapsto T \mapsto S,
\end{equation}
leading to a matrix valued function $S$ that satisfies
\begin{equation} \label{asympS}
    S(z) = I + O\left(\frac{1}{n}\right)
\end{equation}
uniformly for $z \in \mathbb C \setminus \Sigma_S$, where
$\Sigma_S$ is a system of contours shown in
\cite[Figure 13]{KVW}, see \cite[Theorem 6.4]{KVW}.

The starting point of the asymptotic analysis is the Riemann-Hilbert problem for $Y$ on
the contour $\Gamma$ as given in the beginning of Section 2.
We choose $\Gamma$ so that it passes through $\Gamma_P$
and $\Gamma_R$, and the remaining parts are in $D_{\infty,U}$
and $D_{\infty,L}$, see Figure \ref{fig:traject}.
The transformation $Y \mapsto U$ is given by
\begin{equation} \label{fromYtoU}
    U(z) = L^{-n-1} Y(z) \diag\left(e^{-(n+1)g_P(z)}, e^{(n+1)g_A(z)}, e^{-(n+1)g_R(z)}
    \right) L^{n+1}
\end{equation}
where $L = \diag(e^{\ell/3}, e^{-2\ell/3}, e^{\ell/3})$
with $\ell = \log 2 - \pi i$, see \cite[Section 5.2]{KVW}.
This transformation has the effect of normalizing the Riemann-Hilbert
problem at infinity (i.e., $U(z) \to I$ as $z \to \infty$)
and to introduce oscillatory jumps on the contours $\Gamma_P$
and $\Gamma_R$.

The transformation $U \mapsto T$ consists of opening of lenses
around $\Gamma_P$ and $\Gamma_R$. We will not give all the details
about this transformation here, but we refer to \cite[Section 5.3]{KVW}.
The result of the transformation is that $T$ is the solution of a
Riemann-Hilbert problem on a complicated system of contours, but
the jump matrices on all contours tend to the identity matrix
as $n \to \infty$, with the exception of the jump matrices on
$\Gamma_P$ and $\Gamma_R$, which are independent of $n$.

For the final transformation $T \mapsto S$ we surround each of
the branch points $z_j$ by a small disk, and we define
$S = T (M^{(j)})^{-1}$ inside the disk around $z_j$,
and $S= TN^{-1}$ outside of these four disks. The matrix-valued
functions $N$ and $M^{(j)}$, $j=1,2,3,4$ are explicitly constructed
parametrices that resemble $T$ as $n$ gets large.
The result is that $S = I + {\O}(1/n)$ uniformly

Since $X = Y^{-t}$, we get from (\ref{transformationsY})
the corresponding transformations for $X$
\begin{equation} \label{transformationsX}
    X \mapsto U^{-t} \mapsto T^{-t} \mapsto S^{-t},
\end{equation}
and from (\ref{asympS}) we also get $S^{-t} = I + {\O}(1/n)$
uniformly on $\mathbb C \setminus \Sigma_S$.
The proofs of the theorems then follow in the same way
as the proofs in \cite[Section 7]{KVW}.
The role of the exponential prefactors  is now played
by $\widetilde{N} := N^{-t}$, and so to get the explicit
formulas of Theorem 3.4 we need to evaluate the entries
of $\widetilde{N}$.

\subsection{Calculation of $\widetilde{N}$}

We have that
$\widetilde{N} := N^{-t}$ solves the following Riemann-Hilbert
problem  (compare \cite[\S 6.1]{KVW}),
\begin{enumerate}
\item[(1)] $\widetilde{N}$ is defined and analytic in $\mathbb{C} \setminus (\Gamma_P
\cup \Gamma_R)$.
\item[(2)] $\widetilde{N}$ has jumps on $\Gamma_P$ and $\Gamma_R$ given by
\begin{equation}   \label{eq:NjumpP}
   \widetilde{N}_+(z) = \widetilde{N}_-(z) \begin{pmatrix}
                    0 & z^{-1}e^{-3z} & 0 \\
                    -ze^{3z} & 0 & 0 \\
                    0 & 0 & 1
                    \end{pmatrix}, \qquad z \in \Gamma_P,
\end{equation}
and
\begin{equation}   \label{eq:NjumpR}
   \widetilde{N}_+(z) = \widetilde{N}_-(z) \begin{pmatrix}
                    1 & 0 & 0 \\
                    0 & 0 & -ze^{-3z} \\
                    0 & z^{-1}e^{3z} & 0
                    \end{pmatrix}, \qquad z \in \Gamma_R.
\end{equation}
\item[(3)] $\widetilde{N}(z) = I + \O\left(\frac1z\right)$ as $z \to \infty$.
\end{enumerate}

\begin{proposition}  \label{lem:N}
A solution of the Riemann--Hilbert problem for $\widetilde{N}$ is given by
\begin{equation}  \label{eq:NF}
  \widetilde{N}(z) = \begin{pmatrix}
          \widetilde{F}_1(\psi_P(z)) & \widetilde{F}_1(\psi_Q(z)) &
\widetilde{F}_1(\psi_R(z)) \\
          \widetilde{F}_2(\psi_P(z)) & \widetilde{F}_2(\psi_Q(z)) &
\widetilde{F}_2(\psi_R(z)) \\
          \widetilde{F}_3(\psi_P(z)) & \widetilde{F}_3(\psi_Q(z)) &
\widetilde{F}_3(\psi_R(z))
          \end{pmatrix},
\end{equation}
where
\begin{equation}
   \widetilde{F}_1(w) = \frac{-w(w-1)\widetilde{G}(w)}{\sqrt{3 w^4 + 1}},
   \quad
   \widetilde{F}_2(w) = \frac{(w^2-1)\widetilde{G}(w)}{3\sqrt{3 w^4 + 1}},
   \quad \label{eq:F2}
   \widetilde{F}_3(w) = \frac{w(w+1)\widetilde{G}(w)}{\sqrt{3 w^4 + 1}},
\end{equation}
with $\sqrt{3w^4 + 1}$ defined and analytic in
$\mathbb C \setminus (\psi_{P+}(\Gamma_P) \cup \psi_{R+}(\Gamma_R))$, and such that
it is positive for large positive $w$.
The function $\widetilde{G}$ is defined by
\begin{equation} \label{eq:G}
    \widetilde{G}(w)= \left\{ \begin{array}{ll}
    \frac{1}{w}e^{\frac{(w+1)(2w-1)}{w(w-1)}} & \mbox{ for } w\in\psi(\mathcal{R}_P), \\[10pt]
    \left(\frac{w^2-1}{w^2-1/3}\right)e^{\frac{2w^2}{w^2-1}} & \mbox{ for } w \in \psi(\mathcal{R}_Q), \\[10pt]
    \frac{1}{w}e^{\frac{(w-1)(2w+1)}{w(w+1)}} & \mbox{ for } w\in\psi(\mathcal{R}_R).
    \end{array} \right.
\end{equation}
\end{proposition}
\begin{proof}
Rather than taking the inverse transpose of the solution found in
\cite[Proposition 6.1]{KVW}, we prefer to repeat the argument given
there in the present situation. One checks easily that the only
difference lies in replacing the jumps in \cite[Equation (6.14)]{KVW} by
their inverse. Hence the function $G$ used in the proof of \cite[Proposition 6.1]{KVW}
should be replaced with $\widetilde{G}=G^{-1}$.
This leads to the asserted solution for $\widetilde{N}$.
\end{proof}

\begin{corollary}
The entries in the second row of $\widetilde{N}$ can also be
expressed as
\begin{align} \label{eq:tildeN21}
\widetilde{N}_{21}(z) & = \frac{(\psi_P^2(z)-1)^2e^{g_P(z)}}
{2\sqrt{3 \psi_P^4(z) + 1}}, \\
\widetilde{N}_{22}(z) & = - \frac{z^2 (\psi_Q^2(z)-1)^2e^{-g_A(z)}}
{\sqrt{3 \psi_Q^4(z) + 1}}, \label{eq:tildeN22} \\
\widetilde{N}_{23}(z) & = \frac{(\psi_R^2(z)-1)^2e^{g_R(z)}}
{2\sqrt{3 \psi_R^4(z) + 1}}. \label{eq:tildeN23}
\end{align}
\end{corollary}
\begin{proof}
From (\ref{eq:F2}) and the definitions of the functions $F_{2}$
and $G$ in \cite[Proposition 6.1]{KVW}, we derive that
\[
\widetilde{F}_2(w) = \frac{(w^2-1)}{3\sqrt{3 w^4 + 1}G (w)}=
\frac{(w^2-1)^{2}}{( 3 w^4 + 1)F_{2}(w)}.
\]
Now, replacing $w$ with $\psi_{P} (z)$, $\psi_{Q}
(z)$, $\psi_{R} (z)$, respectively, and using the expressions of $N_{21}
(z)=F_{2} (\psi_{P} (z))$, $N_{22} (z)=F_{2} (\psi_{Q} (z))$,
$N_{23} (z)=F_{2} (\psi_{R} (z))$ that are given in \cite[Equation
(6.20)]{KVW}, and also noting that $g_A = g_P + g_R$, we arrive at
(\ref{eq:tildeN21})--(\ref{eq:tildeN23}).
\end{proof}

\subsection{Proof of Theorem \ref{asymp}}
We get from (\ref{fromYtoU}) and (\ref{asympS}), and the fact
that $U = T = SN$ away from $\Gamma_P$ and $\Gamma_R$, that
\begin{equation} \label{asympX1}
    X = Y^{-t} =  L^{-n-1} \left(I + {\O}(1/n)\right)  \widetilde{N}  \diag\left(e^{-(n+1)g_P}, e^{(n+1)g_A}, e^{-(n+1)g_R}
    \right) L^{n+1}
\end{equation}
uniformly in compact subsets of $\mathbb C \setminus (\Gamma_P \cup \Gamma_R)$.

From (\ref{asympX1}) and (\ref{eq:Xout})-(\ref{eq:Xin}) we get that
\[ A_n(z) = X_{22}(z) = \widetilde{N}_{22}(z) e^{(n+1) g_A(z)} \left(1 + {\O}\left(\frac{1}{n}\right)\right), \]
uniformly for $z$ in compact subsets of $\mathbb C \setminus (\Gamma_P \cup \Gamma_R)$,
which by (\ref{eq:tildeN22}) leads to (\ref{asympA}).

In the same way, we get from (\ref{asympX1}) and (\ref{eq:Xout}) that
\[ B_n(z) = X_{21}(z) z^{3n+2} = \widetilde{N}_{21}(z) (-2)^{n+1} z^{3n+2} e^{-(n+1)g_P(z)}
    \left(1+{\O}\left(\frac{1}{n}\right)\right), \]
uniformly for $z$ in compact subsets of the exterior of $\Gamma$.
Using (\ref{eq:tildeN21}) and the fact that $g_B(z) + g_P(z) = 3\log z$ for
$z \in \mathbb C \setminus D_P^*$, we find the formula (\ref{asympB}) for $z$
in the exterior of $\Gamma$. For $z$ inside $\Gamma$, we have
by (\ref{asympX1}) and (\ref{eq:Xin})
\begin{equation}  \begin{aligned}
    B_n(z) & = A_n(z) e^{-3nz} - E_n^{(1)}(z)
    = X_{22}(z) e^{-3nz} + X_{21}(z) z^{3n+2} \\
    & = \widetilde{N}_{22}(z) e^{(n+1)g_A(z)} e^{-3nz} \left(1 + {\O}\left(\frac{1}{n}\right)\right) \\
     & \qquad + \widetilde{N}_{21}(z) (-2)^{n+1} z^{3n+2} e^{-(n+1)g_P(z)}
    \left(1+{\O}\left(\frac{1}{n}\right)\right).
    \label{asympB2}
    \end{aligned}
    \end{equation}
The first term in (\ref{asympB2}) dominates for $z \in D_P^*$,
and this leads to the second asymptotic formula in (\ref{asympB}), since
$g_A(z) - 3z = g_B + \ell$ for $z \in D_P^*$. The second term in (\ref{asympB2})
dominates in the remaining part inside $\Gamma$ and this again leads to
the first asymptotic formula in (\ref{asympB}). A more detailed analysis
for $z$ near $\Gamma_P$ and $\Gamma_R$ reveals that (\ref{asympB}) is
uniformly valid near $\Gamma_P$ and $\Gamma_R$ as well.

The proofs of (\ref{asympC}), (\ref{asympE1}), and (\ref{asympE2})
are similar.
\qed

\subsection{Proof of Theorem \ref{lim-zeros}}
The limits for the counting measures $\nu_{A_{n}}$, $\nu_{B_{n}}$, and
$\nu_{C_{n}}$ follow directly from the strong convergence results in Theorem \ref{asymp} and the unicity
theorem for logarithmic potentials, see e.g.\
\cite[Theorem II.2.1]{ST}.

To establish the limits for $\nu_{E_{n}^{(1)}}$ and $\nu_{E_{n}^{(2)}}$
we follow the proof of \cite[Theorem 2.5]{KVW} given in Section 7.2 of
\cite{KVW}, which established the limit for $\nu_{E_{n}}$ as $n \to \infty$.
In following this proof we have to make some obvious modifications
to see that the proof applies to the present situation
as well.
The only modification that may not be immediately obvious is
that in \cite[Equation (7.13)]{KVW} we used an integral representation for $E_n$ to conclude that
$\frac{1}{n} \left.\log(z^{-3n-2} E_n(z)) \right|_{z=0}$
is bounded from below as $n \to \infty$. We do not have
an integral representation for $E_n^{(1)}$ and $E_n^{(2)}$,
but the fact that
$\frac{1}{n} \left.\log(z^{-3n-2} E_n^{(j)}(z)) \right|_{z=0}$
is bounded from below as $n \to \infty$ for $j=1,2$, follows easily
from the asymptotic formulas (\ref{asympE1}) and (\ref{asympE2}).
\qed

\subsection{Proof of Theorem \ref{asym-curves}}
This is similar to the proofs of Theorems 2.9 and 2.10 in
\cite[Sections 7.3 and 7.4]{KVW}.
We omit the details.
\qed

\subsection{Proof of Theorem \ref{asym-branch}}
The proof is similar to the proof of Theorem 2.11 in \cite[Section 7.5]{KVW}.
Applying the inverse transform to \cite[Equation (7.19)]{KVW}
leads to an expression of $X(z)$ involving the matrices $\widetilde{N}$ and
$(S^{-1})^t=I+{\O}(1/n)$. Multiplying that expression by the unit vector
$\begin{pmatrix} 1 & 0 & 0 \end{pmatrix}^t$ as in the proof of \cite[Theorem
2.11]{KVW} leads after some calculations to the asymptotic formula (\ref{asympbrE1}) for
$E_n^{(1)}$, while multiplying by the unit vector
$\begin{pmatrix} 0 & 1 & 0 \end{pmatrix}^t$ leads to the formula (\ref{asympbrA})
for $A_n$. In order to get the asymptotics (\ref{asympbrB})
for $B_n$, one has to choose $z$ outside of the contour $\Gamma$ for the
Riemann-Hilbert problem for $X$. Then one performs similar calculations
in the regions $f_1^{-1}(I)$ or $f^{-1}(II)$, as they were defined while
constructing the parametrices near the branch points in \cite[\S 6.2]{KVW},
see also \cite[Figure 12]{KVW}. We note that the following relation between
the Airy function and its derivative is used in the above computations,
$$\Ai(s)\Ai'(\omega s)-\omega^2 \Ai(\omega s)\Ai'(s)=-\frac{e^{i\pi/6}}{2\pi},$$
where $\omega=e^{2\pi i/3}$, see e.g. \cite[Ex. 8.1, p. 416]{OLV}.
\qed

\subsection{Proof of Corollary \ref{asym-zero}}
The formulas for the zeros of $A_n$, $B_n$, and $E_n^{(1)}$ near
$z_1$ follow from the asymptotics near the branch points given in
Theorem \ref{asym-branch}. We refer to the proof of
\cite[Corollary 2.12]{KVW} given in \cite[Section 7.6]{KVW}
for details.
\qed

\end{document}